\documentclass[]{amsart}

\usepackage{amssymb}
\usepackage{amsmath}

\usepackage[all]{xy}

\newcommand{\Norm}{{\mathcal N}}
\newcommand{\F}{{\mathcal F}}
\newcommand{\Ol}{{\mathcal O}}
\newcommand{\f}{\varphi}
\newcommand{\B}{{\mathcal B}}
\newcommand{\E}{{\mathcal E}}
\newcommand{\Di}{{\overline{D}}}
\newcommand{\V}{{\mathcal V}}

\newcommand{\acca}{\widetilde{H}}
\newcommand{\pu}{{\mathbb P^1}}
\newcommand{\proj}{\mathbb P}
\newcommand{\quadr}{\mathbb Q}
\newcommand{\pt}{{\mathbb P^3}}
\newcommand{\pd}{{\mathbb P^2}}
\newcommand{\pdue}{\mathbb P^2}
\DeclareMathOperator{\loc}{\mathrm{Locus}}
\DeclareMathOperator{\cloc}{\mathrm{ChLocus}}
\newcommand{\ratcurves}{\textrm{Ratcurves}^n(X)}
\newcommand{\ratcurvesx}{\textrm{Ratcurves}^n(X,x)}
\newcommand{\om}{\textrm{Hom}}
\newcommand{\Aut}{\textrm{Aut}}
\newcommand{\Univ}{\textrm{Univ}}
\newcommand{\cone}{\textrm{NE}}
\DeclareMathOperator{\pic}{Pic}
\newcommand{\sa}{\hspace{1pt}}
\newcommand{\W}{{\overline{W}}\sa}
\newcommand{\onespan}[1]{\langle #1 \rangle}
\newcommand{\twospan}[2]{\langle #1,#2 \rangle}
\newcommand{\threespan}[3]{\langle #1,#2,#3 \rangle}
\newcommand{\morespan}[2]{\langle #1, \ldots, #2 \rangle}

\newtheorem{theorem}{Theorem}[section]
\newtheorem{lemma}{Lemma}[section]
\newtheorem{proposition}{Proposition}[section]
\newtheorem{corollary}{Corollary}[section]

\theoremstyle{definition}
\newtheorem{definition}{Definition}[section]
\newtheorem{example}{Example}[section]

\theoremstyle{remark}
\newtheorem{remark}{Remark}[section]

\allowdisplaybreaks                          
\makeatletter
\@addtoreset{equation}{section}
\makeatother

\begin{document}
\markboth{Elena Chierici and Gianluca Occhetta}
{The cone of curves of Fano varieties of coindex four}

\title{The cone of curves of Fano varieties of coindex four}

\author{Elena Chierici and Gianluca Occhetta}

\address{Dipartimento di Matematica, via Sommarive 14, I-38050 Povo (TN)}
\email{e.chierici@email.it, gianluca.occhetta@unitn.it} 

\begin{abstract}
We classify the cones of curves of Fano varieties of dimension greater or equal
than five and (pseudo)index $\dim X -3$, describing the number and type of their
extremal rays.

\end{abstract}
\maketitle

\section{Introduction} 

A smooth complex projective variety is called {\sf Fano} if its
anticanonical bundle $-K_X$ is ample.\\ 
The {\sf index} of $X$, $r_X$, is the largest natural
number $m$ such that $-K_X=mH$ for some (ample) divisor $H$ on $X$, while the 
{\sf coindex} of $X$ is defined as $\dim X + 1 -r_X$.\\
Since $X$ is smooth, $\pic(X)$ is torsion free and therefore the divisor
$L$ satisfying $-K_X = r_XL$ is uniquely determined and called the 
{\sf fundamental divisor} of $X$.\\
It is known that $0 < r_X \le \dim X+1$ and, by a theorem of 
Kobayashi and Ochiai \cite{Koba}, $r_X= \dim X +1$ if and only if 
$(X,L) \simeq (\proj^{\dim X}, \Ol_{\proj}(1))$, and $r(X)= \dim X$ if and only if 
$(X,L) \simeq (\quadr^{\dim X}, \Ol_{\quadr}(1))$.\\
Fano varieties of coindex two, which are called {\sf del Pezzo} varieties,
have been classified in \cite{Fu} using the Apollonius method, i.e. proving that
the linear sistem $|L|$ contains a smooth divisor and constructing a ladder
down to the well-known case of surfaces.\\
The same method works for Fano varieties of coindex 
three, called {\sf Mukai} varieties; in \cite{Mu} Mukai announced
the classification assuming the existence of a smooth member in $|L|$,
and this was proved by Mella in \cite{Me}.\par

\smallskip
Since the classification of Fano fourfolds is very far from being known,
it is not possible to use Apollonius method to study Fano varieties
of coindex four; however, with different techniques which involve
the study of families of rational curves, it is possible to describe their 
structure, i.e. their cone of curves.\\
Families of rational curves are related to another invariant of a Fano variety,
the {\sf pseudoindex} $i_X$, introduced by Wi\'sniewski in \cite{Wimu}, which is
defined as the minimum anticanonical degree of rational curves on $X$.\\
The pseudoindex is related to the dimension and to the Picard number $\rho_X$
of a Fano variety by a conjecture of Mukai \cite{Kata}, which states that
\begin{equation} \label{conj}
\rho_X(i_X -1) \le \dim X,
\end{equation}
equality holding if and only if $X \simeq (\proj^{i_X-1})^{\rho_X}$.\\
This conjecture has been recently proved in \cite{ACO} for Fano fivefolds
and for Fano varieties of pseudoindex $i_X \ge \frac{\dim X +3}{3}$ which 
admit a covering unsplit family of rational curves
(this is always the case if $i_X=\dim X-3$ and $\dim X \ge 6$),
so it provides a good bound on the Picard number
of Fano varieties of pseudoindex $i_X \ge \dim X -3$ and dimension
greater than four.\\

By definition, the pseudoindex is an integral multiple of $r_X$, so varieties of coindex four have
pseudoindex $i_X \ge \dim X-3$; by the main result in \cite{AOlong} it is easy to
prove that, for varieties of dimension $\ge 5$, Picard number $\ge 2$ and coindex four,
index and pseudoindex coincide. For this reason the classification of the cone of curves of Fano varieties
of coindex four and dimension $\ge 5$ can be regarded as a special case of the following

\begin{theorem}\label{main}
Let $X$ be a Fano variety of dimension $n \ge 5$, pseudoindex $i_X = n -3$ and Picard
number $\rho_X \ge 2$. Then $\cone(X)$ is generated by $\rho_X$ rays.\\ 
More precisely, we have the following list of possibilities, 
where $F$ stands for a fiber type contraction, $D_i$ for a divisorial contraction whose
exceptional locus is mapped to a $i-$dimensional subvariety and $S$ for a small
contraction.\\
All cases are effective.
\end{theorem}

\begin{center}
\begin{tabular}{|c|c||c|c|c|c|c|}\hline
\quad $\dim X$ \quad & \quad $\rho_X$ \quad &  \quad $R_1$
\quad & \quad $R_2$ \quad & \quad $R_3$ \quad & \quad $R_4$ \quad & \quad $R_5$
\quad \\
\hline\hline
 $5 $& $2$  & $F$ & $F$  & &  & \\
   &   & $F$ & $D_0$  & & & \\
  &    & $F$ & $D_1$ & & & \\
  
 &     & $F$ & $D_2$  & & & \\
 &     & $F$ & $S$  & & & \\
	&  & $D_2$ & $D_2$ & & & \\
	&  & $D_2$ & $S$ & & & \\
\hline
& $3$ & $F$ & $F$  & $F$ & & \\
&      & $F$ & $F$ & $S$ & & \\
 &    & $F$ & $F$ & $D_1$ & & \\
 &    & $F$ & $F$ & $D_2$ & & \\
 &     & $F$ & $D_2$ & $D_2$ & & \\
\hline
 & $4$  & $F$ & $F$  & $F$ & $F$ &\\
 &     & $F$ & $F$ & $F$ & $D_2$ & \\
\hline
 & $5$  & $F$ & $F$  & $F$ & $F$ & $F$ \\
 \hline\hline
 $6$ & $2$  & $F$ & $F$  &  &&  \\
     &   & $F$ & $D_1$  & &&   \\
     &   & $F$ & $D_2$  & &&   \\
     &   & $F$ & $S$   & &&   \\
	 \hline
	 & $3$  & $F$ & $F$ & $F$ &&  \\
\hline\hline
 $7$ & $2$  & $F$ & $F$ & &&  \\
 	 &   & $F$ & $D_2$ & &&   \\
\hline\hline
 $8$ & $2$  & $F$ & $F$ & &&  \\
\hline
\end{tabular}
\end{center}~\par

\medskip
Note that, under the stronger assumption that $r_X=\dim X - 3$, a local description 
of the Fano-Mori contractions of $X$ has been achieved by Andreatta and Wi{\'s}niewski 
in \cite{AW}, hence we expect that in this case it will be possible to
go further in the direction of an effective classification. 

\medskip
The paper is organized as follows: in section two we collect basic material
concerning Fano-Mori contractions and families of rational curves as well as 
some definitions and results in \cite{ACO} which we will use extensively
throughout the proof, while in section three we construct examples showing 
that all cases in our list are effective.\\
We then start the proof of theorem \ref{main}: the case $\dim X=5$, which is the hardest, 
is divided into two parts: in section four we deal with Fano fivefolds which admit a quasi-unsplit 
locally unsplit covering family, which turn out to have always a fiber type contraction,
and in section five we study the remaining cases, proving the following\par 

\smallskip
\begin{theorem}\label{shelob} 
Let $X$ be a Fano fivefold of pseudoindex two which does not have a 
covering quasi-unsplit locally unsplit family of rational curves; then $\rho_X = 2$ and $X$ 
is  the blow-up of $\proj^5$ along a two-dimensional smooth quadric (a section of $\Ol(2)$
in a linear $\proj^3 \subset \proj^5$), 
or along a cubic scroll ($\proj_\pu(\Ol(1) \oplus \Ol(2))$ embedded in an hyperplane of $\proj^5$
by the tautological bundle), or along a Veronese surface.
\end{theorem}

Finally, the last section contains the proof of theorem \ref{main} in the easier case
of varieties of dimension greater than five.


\section{Background material}

\subsection{Extremal contractions} 
Let $X$ be a smooth complex projective variety
of dimension $n$ and let $K_X$ be its canonical divisor.\\
By Mori's {\sf Cone Theorem} the closure of the cone of effective 1-cycles into
the $\mathbb R$-vector space of 1-cyles modulo numerical equivalence,
$\overline {\cone(X)} \subset N_1(X)$, is locally polyhedral in the part contained in
the set $\{Z\in N_1(X) \ |\ K_X \cdot Z<0\}$; an {\sf extremal face} $\sigma$ 
of $X$ is a face of this locally polyhedral part and
an {\sf extremal ray} is an extremal face of dimension one.\\
Note that, if $X$ is a Fano variety, then $\overline{\cone(X)}=\cone(X)$ is polyhedral
and any face of $\cone(X)$ is an extremal face.\\ 
To every extremal face one can
associate a morphism to a normal variety; namely we have the following
{\sf Contraction Theorem} due to Kawamata and Shokurov:

\begin{theorem}\label{contractionth}
Let $X$ and $\sigma$ be as above.
Then there exists a projective morphism $\f:X \to W$ from $X$ onto a normal
variety $W$ which is characterized by the following properties:
\begin{itemize}
\item[{i})]  for every irreducible curve $C$ in $X$, $\f(C)$
is a point if and only if the numerical class of $C$ is in $\sigma$;
\item[{ii})] $\f$ has connected fibers.
\end{itemize}
\end{theorem}

\smallskip
\begin{definition}
The map $\f$ of the above theorem is usually called the {\sf Fano-Mori
contraction}
(or the {\sf extremal contraction}) associated to the face $\sigma$.
A Cartier divisor $H$ such that $H = \f ^*A$ for an ample divisor $A$ on $W$
is called a {\sf good supporting divisor} of the map $\f$ (or of the face $\sigma$).\\
An extremal ray $R$ is called {\sf numerically effective}, or of
{\sf fiber type}, if $\dim W < \dim X$, otherwise the ray is {\sf non nef} or {\sf birational};
the terminology is due to the fact that if $R$ is non nef then there exists an irreducible divisor 
$D_R$ which is negative on curves in $R$.\\
We usually denote with $E = E(\f):= \{ x \in X\ |\  \dim(\f^{-1} \f(x)) > 0\}$
the {\sf exceptional locus} of $\f$; if $\f$ is of fiber type then of course $E=X$.\\
If the codimension of the exceptional locus of a birational ray $R$ is equal to one,
the ray and the associated contraction are called {\sf divisorial}, otherwise they are called {\sf small}.
\end{definition}


\bigskip
\subsection{Families of rational curves} 
For this subsection our main reference is \cite{Kob}, 
with which our notation is coherent.\\
Let $X$ be a normal projective variety and let $\om(\pu,X)$ be the scheme parametrizing
morphisms $f: \pu \to X$; we consider the open subscheme $\om_{bir}(\pu,X) \subset \om(\pu,X)$, 
corresponding to those morphisms which are birational onto their image, and its normalization
$\om^n_{bir}(\pu,X)$; the group $\Aut(\pu)$ acts on $\om^n_{bir}(\pu,X)$ and the
quotient exists.\par

\smallskip
\begin{definition}
The space $\ratcurves$ is the quotient of $\om^n_{bir}(\pu,X)$ by $\Aut(\pu)$, and the space
$\Univ(X)$ is the quotient of the product action of $\Aut(\pu)$ on
$\om^n_{bir}(\pu, X)\times\pu$.
\end{definition}

\smallskip
\begin{definition} \label{Rf}
We define a {\sf family of rational curves} to be an irreducible component
$V \subset \ratcurves$.\\
Given a rational curve $f:\pu \to X$ we will call a {\sf family of
deformations} of $f$ any irreducible component $V \subset
\ratcurves$ containing the equivalence class of $f$.
\end{definition}

\medskip
Given a family $V$ of rational curves, we have the following basic diagram:
$$
\xymatrix{p^{-1}(V)=:U  \ar[r]^(.65){i} \ar[d]_{p} & X\\
 V & & }$$
where $i$ is the map induced by the evaluation $ev:\om^n_{bir} (\pu, X)\times \pu \to X$
and $p$ is a $\pu$-bundle. \\
We define $\loc(V)$ to be the image of $U$ in $X$;
we say that $V$ is a {\sf covering family} if $\overline{\loc(V)}=X$.
We will denote by $\deg V$ the anticanonical degree of the family $V$, i.e. the integer
$-K_X \cdot C$ for any curve $C \in V$.\par
\smallskip
Given a family $V \subseteq \ratcurves$ and a point $x \in \loc(V)$. we denote by $V_x$ the subscheme
of $V$ parametrizing rational curves passing through $x$.\par

\begin{definition} 
Let $V$ be a family of rational curves on $X$. Then
\begin{itemize}
\item[(a)] $V$ is {\sf unsplit} if it is proper; 
\item[(b)] $V$ is {\sf locally unsplit} if for the general $x \in \loc(V)$ every
component $V_x$ of \linebreak
$V \cap \ratcurvesx$ is proper;
\item[(c)] $V$ is {\sf generically unsplit} if there is at most a finite number of curves of $V$
passing through two general points of $\loc(V)$.
\end{itemize}
\end{definition}

\smallskip
\begin{remark} 
Note that $(a) \Rightarrow (b) \Rightarrow (c)$.
\end{remark}

\smallskip
\begin{example}\label{unex} 
Let $R_i=\mathbb R_+[C_i]$ be an extremal ray such
that the anticanonical degree of $[C_i]$ is minimal in $R_i$; $C_i$ is often called a 
{\sf minimal extremal rational curve}.\\
If we denote by $R^i$ an irreducible component of $\ratcurves$ containing $C_i$, then
the family $R^i$ is unsplit: in fact, if $C_i$ degenerates into a reducible cycle, its components
must belong to the ray $R_i$, since $R_i$ is extremal; but in $R_i$
the curve $C_i$ has the minimal intersection with the anticanonical
bundle, hence this is impossible. 
\end{example}

\smallskip
\begin{proposition}(IV.2.6 in \cite{Kob})\label{iowifam} 
Let $X$ be a smooth projective variety
and $V$ a family of rational curves.\\
Assume either that $V$ is generically unsplit and $x$ is a general point in $\loc(V)$ or that
$V$ is unsplit and $x$ is any point in $\loc(V)$. Then
 \begin{itemize}
      \item[(a)] $\dim X + \deg V \le \dim \loc(V)+\dim \loc(V_x) +1$;
      \item[(b)] $\deg V \le \dim \loc(V_x)+1$.
   \end{itemize}
\end{proposition}

This last proposition, in case $V$ is the unsplit family of deformations of a minimal extremal
rational curve, gives the {\sf fiber locus inequality}:

\begin{proposition}\label{fiberlocus} Let $\f$ be a Fano-Mori contraction
of $X$ and let $E = E(\f)$ be its exceptional locus;
let $S$ be an irreducible component of a (non trivial) fiber of $\f$. Then
$$\dim E + \dim S \geq \dim X + l -1$$
where 
$$l =  \min \{ -K_X \cdot C\ |\  C \textrm{~is a rational curve in~} S\}.$$
If $\f$ is the contraction of a ray $R$, then $l$ is called the {\sf length of the ray}.
\end{proposition}

\smallskip
\begin{definition}\label{CF}
We define a {\sf Chow family of rational curves} to be an irreducible component
$\V \subset \textrm{Chow}(X)$ parametrizing rational and connected 1-cycles.
\end{definition}

\medskip
Given a Chow family of rational curves, we have a diagram as before, coming from the universal
family over $\textrm{Chow}(X)$.
\begin{equation}\label{chowdiagram}
\xymatrix@=30pt{\mathcal U  \ar[r]^{i} \ar[d]_{p} & X\\
 \V & & }
\end{equation}
In the diagram $i$ is the map induced by the evaluation and the fibers of $p$ are connected and have
rational components. Both $i$ and $p$ are proper (see for instance II.2.2 in \cite{Kob}).

\smallskip
\begin{definition}
If $V$ is a family of rational curves, the closure of the image of
$V$ in $\textrm{Chow}(X)$ is called the {\sf Chow family associated to} $V$.
\end{definition}

\smallskip
\begin{remark}
If $V$ is proper, i.e. if the family is unsplit, then $V$ corresponds to the normalization
of the associated Chow family $\V$.
\end{remark}

\smallskip
\begin{definition}
Let $\V$ be the Chow family associated to a family of rational curves $V$. We say that
$V$ is {\sf quasi-unsplit} if every component of any reducible cycle in $\V$ is
numerically proportional to $V$.
\end{definition}

\bigskip
\subsection{Notation}

From now on we will denote with\\

$V$ a family of rational curves;

$\V$ the associated Chow family;

$[V]$ the numerical equivalence class in $N_1(X)$ of a general curve belonging to the family $V$;

$R_i$ an extremal ray of $X$;

$R^i$ the (unsplit) family of deformations of a minimal rational curve in $R_i$ (see example \ref{unex});

$\f_{R_i}$ or $\f_i$ the extremal contraction associated with the ray $R_i$.

\bigskip
\subsection{Chains of rational curves} 

For all the missing proofs in the rest of the section we refer the reader to \cite{ACO}.\par

\medskip
Let $X$ be a smooth variety, $\V^1, \dots, \V^k$ Chow families of rational curves
on $X$ and $Y$ a subset of $X$.\par

\smallskip
\begin{definition}
We denote by $\loc(\V^1, \dots, \V^k)_Y$ the set of points $x \in X$ such that there exist
cycles $C_1, \dots, C_k$ with the following properties:
   \begin{itemize}
      \item $C_i$ belongs to the family $\V^i$;
      \item $C_i \cap C_{i+1} \not = \emptyset$;
      \item $C_1 \cap Y \not = \emptyset$ and $x \in C_k$,
   \end{itemize}
i.e. $\loc(\V^1, \dots, \V^k)_Y$ is the set of points that can be joined to $Y$ by a connected
chain of $k$ cycles belonging \underline{respectively} to the families $\V^1, \dots, \V^k$.
\end{definition}

Note that $\loc(\V^1, \dots, \V^k)_Y \subset \loc(\V^k)$.\par

\smallskip
\begin{remark} \label{Lr?}
If $Y$ is a closed subset, then $\loc(\V^1, \dots, \V^k)_Y$ is
closed.
\end{remark}

\smallskip

\begin{definition} Let $V^1, \dots, V^k$ be unsplit families on $X$.
We will say that $V^1, \dots, V^k$ are {\sf numerically independent} if their numerical classes
$[V^1], \dots ,[V^k]$ are linearly independent in the vector space $N_1(X)$.
If moreover $C \subset X$ is a curve we will say that $V^1, \dots, V^k$ are numerically 
independent from $C$ if the class of $C$ in $N_1(X)$ is not contained in the
vector subspace generated by $[V^1], \dots ,[V^k]$.   
\end{definition}

\begin{lemma} \label{locy}
Let $Y \subset X$ be a closed subset and $V$ an unsplit family.
Assume that curves contained in $Y$ are numerically independent from curves in $V$, and that
$Y \cap \loc(V) \not= \emptyset$. Then for a general $y \in Y \cap \loc(V)$
\begin{itemize}
      \item[(a)] $\dim \loc(V)_Y \ge \dim (Y \cap \loc(V)) + \dim \loc(V_y);$
      \item[(b)] $\dim \loc(V)_Y \ge \dim Y + \deg V - 1$.
\end{itemize}
Moreover, if $V^1, \dots, V^k$ are numerically independent
unsplit families such that curves contained in $Y$ are numerically independent
from curves in $V^1, \dots, V^k$ then either\linebreak
$\loc(V^1, \ldots, V^k)_Y=\emptyset$ or
\begin{itemize}
      \item[(c)] $\dim \loc(V^1, \ldots, V^k)_Y \ge \dim Y +\sum \deg V^i -k$.
\end{itemize}
\end{lemma}

\smallskip
\begin{definition}
We denote by $\cloc_m(\V^1, \dots, \V^k)_Y$ the set of points $x \in X$
such that there exist cycles $C_1, \dots, C_m$ with the following properties:
   \begin{itemize}
      \item $C_i$ belongs to a family $\V^j$;
      \item $C_i \cap C_{i+1} \not = \emptyset$;
      \item $C_1 \cap Y \not = \emptyset$ and $x \in C_m$,
   \end{itemize}
i.e. $\cloc_m(\V^1, \dots, \V^k)_Y$ is the set of points that can be joined to $Y$ by a
connected chain  of at most $m$ cycles belonging to the families $\V^1, \dots, \V^k$.
\end{definition}

\smallskip
\begin{definition}
We define a relation of {\sf rational connectedness with respect to $\V^1, \dots, \V^k$}
on $X$ in the following way: $x$ and $y$ are in rc$(\V^1,\dots,\V^k)$-relation if there
exists a chain of rational curves in $\V^1, \dots ,\V^k$ which joins $x$ and $y$, i.e.
if $y \in \cloc_m(\V^1, \dots, \V^k)_x$ for some $m$.
\end{definition}
\smallskip
To the rc$(\V^1,\dots,\V^k)$-relation we can associate a fibration, at least on an open subset.

\begin{theorem}(\cite{Cam}, IV.4.16 in \cite{Kob}) \label{rcvfibration}
There exist an open subvariety $X^0 \subset X$ and a proper morphism with connected fibers
$\pi:X^0 \to Z^0$ such that
   \begin{itemize}
      \item[(a)] the rc$(\V^1,\dots,\V^k)$-relation restricts to an equivalence relation on $X^0$;
      \item[(b)] the fibers of $\pi$ are equivalence classes for the rc$(\V^1,\dots,\V^k)$-relation;
      \item[(c)] for every $z \in Z^0$ any two points in $\pi^{-1}(z)$ can be connected by a chain
      of at most $2^{\dim X - \dim Z}-1$ cycles in $\V^1, \dots, \V^k$.
   \end{itemize}
\end{theorem}

\begin{definition}
In the above assumptions, if $\pi$ is the constant map we say that $X$ is {\sf rc$(\V^1,\dots,\V^k)$-connected}.
\end{definition}

\bigskip
\subsection{Bounding Picard numbers}

In this subsection we list some conditions
under which the numerical class (in $X$) of every curve lying in some subvariety $S \subset X$ 
is contained in a linear subspace of $N_1(X)$ or in a subcone of $\cone(X)$.\\
We write $N_1(S)=\morespan{[V^1]}{[V^k]}$ if the numerical class in $X$ of every curve $C \subset S$ can
be written as $[C]= \sum_i a_i [C_i]$, with $a_i \in \mathbb Q$ and $C_i \in V^i$, and  
$\cone(S)=\morespan{[V^1]}{[V^k]}$ if the numerical class in $X$ of every curve $C \subset S$ can
be written as $[C]= \sum_i a_i [C_i]$, with $a_i \in \mathbb Q_{\ge 0}$ and $C_i \in V^i$.\par

\smallskip
\begin{lemma} (Lemma 1 in \cite{Op}) \label{numequns}
Let $Y \subset X$ be a closed subset and $V$ an unsplit family of rational curves. 
Then every curve contained in $\loc(V)_Y$ is numerically equivalent to
a linear combination with rational coefficients
   $$\lambda C_Y + \mu C_V,$$
where $C_Y$ is a curve in $Y$, $C_V$ belongs to the family $V$ and $\lambda \ge 0$.
\end{lemma}

\smallskip
\begin{corollary}\label{rho1}
Let $V$ be a family of rational curves and $x$ a point in $X$ such that $V_x$ is unsplit.
Then $N_1(\loc(V_x))=\cone(\loc(V_x))=\onespan{[V]}$.
\end{corollary}

\smallskip
\begin{corollary}\label{ray}
Let $R_1$ be an extremal ray of $X$,
$R^1$ a family of deformations of a minimal extremal curve in $R_1$, $x$ a point in $\loc(R^1)$ and
$V$ an unsplit family of rational curves, independent from $R^1$.\\
Then $\cone(\cloc_m(V)_{\loc(R^1_x)})=\twospan{[V]}{R_1}$.
\end{corollary}

\begin{proof} 
Since
   $$\cloc_m(V)_{\loc(R^1_x)} = \loc(V)_{\cloc_{m-1}(V)_{\loc(R^1_x)}},$$
iterating lemma \ref{numequns} $m$ times any curve $C$ in
$\cloc_m(V)_{\loc(R^1_x)}$ can be written as
 $$C \equiv \lambda C_1 + \mu C_V$$
with $\lambda \ge 0$, so we have only to prove that $\mu \ge 0$.\\
If $\mu <0$, then we can write $C_1 \equiv \alpha C_V + \beta
C$ with $\alpha, \beta \ge 0$; but since $C_1$ is extremal this implies that both $[C]$ and $[C_V]$ 
belong to $R_1$, a contradiction.
\end{proof}

\smallskip
\begin{remark}\label{face}
More generally, if $\sigma$ is an extremal face of $\cone(X)$, $F$ is a fiber
of the associated contraction and $V$ is an unsplit family independent from $\sigma$,
the same proof shows that
   $$\cone(\loc(V)_F) = \langle \sigma, [V] \rangle.$$
\end{remark}

\bigskip
\subsection{Rational curves on Fano varieties} \label{fano}

Let $X$ be a Fano variety and $\pi:X^0
\to Z^0$ a proper surjective morphism on a smooth quasiprojective variety $Z^0$ of
positive dimension.\\
By Theorem 2.1 in \cite{KoMiMo} we know that for a general point $z \in Z^0$ there exists a
rational curve $C$ on $X$ of anticanonical degree $\le \dim X +1$
which meets $\pi^{-1}(z)$ without being contained in it (an {\sf horizontal curve}, for short).\\
We consider all the families containing these horizontal curves and,
since they are only a finite number, we have
that the locus of at least one of them dominates $Z^0$.\par

\smallskip
\begin{definition} 
A {\sf minimal horizontal dominating family}
for $\pi$ is a family $V$ of horizontal curves such that $\loc(V)$
dominates $Z^0$ and $\deg V$ is minimal among the families
with this property.\\
If $\pi$ is the identity map we say that $V$ is a {\sf minimal covering family}.
\end{definition}

\begin{lemma} \label{horizontal}
Let $X$ be a Fano variety, and let $\pi:\xymatrix@1{X~ \ar@{-->}[r]& ~Z}$ be the
rationally connected fibration associated to $k$ Chow families $\V^1, \ldots, \V^k$; 
let $V$ be a minimal horizontal dominating family for $\pi$. Then
   \begin{itemize}
      \item[(a)] curves parametrized by $V$ are numerically independent from curves contracted
      by $\pi$;
      \item[(b)] $V$ is locally unsplit;
      \item[(c)] if $x$ is a general point in $\loc(V)$ and $F$ is the fiber containing $x$, then
         $$\dim (F \cap \loc(V_x))=0.$$
   \end{itemize}
\end{lemma}

\smallskip
\begin{remark}\label{rcvdim} 
Let $X$ be a Fano variety, $V^1, \dots, V^k$ locally unsplit families
of rational curves such that $V^1$ is covering and $V^i$ is
horizontal and dominating with respect to the rc$(\V^1, \dots,
\V^{i-1})$-fibration; let $\pi:\xymatrix@1{X~ \ar@{-->}[r]& ~Z}$ be the
rc$(\V^1, \dots, \V^k)$-fibration. Then for general $x_i \in \loc(V^i)$
\begin{equation*}
\sum \dim \loc(V^i_{x_i}) \le \dim X - \dim Z.
\end{equation*}
\end{remark}

\smallskip
\begin{lemma}\label{codim}
Let $X$ be a Fano variety, $V^1, \dots, V^k$  locally unsplit
families of rational curves such that $V^1$ is covering and $V^i$
is horizontal and dominating with respect
to the rc$(\V^1, \dots, \V^{i-1})$-fibration.\\
Let $\pi:\xymatrix@1{X~ \ar@{-->}[r]& ~Z}$
 be the rc$(\V^1, \dots, \V^k)$-fibration and suppose that $\dim Z >0$.\\
Then either $[V^1], \dots, [V^k]$ are contained in an extremal face of
$\cone(X)$ or there exists a small extremal ray $R$ whose exceptional locus is
contained in the indeterminacy locus of $\pi$.
\end{lemma}

\begin{proof}
Since $X$ is normal and $Z$ is proper, the indeterminacy locus $E$ of $\pi$ in $X$ has codimension
$\ge 2$ (see [1.39] in \cite{De}). Take a very ample divisor $H$ on $Z$ and pull
it back to $X$: then  $\pi^*H$ is zero on curves in $\morespan{[V^1]}{[V^k]}$ 
and it is non negative outside the indeterminacy locus of $\pi$.\\
Therefore, either $\pi^*H$ is nef on $X$ and $\morespan{[V^1]}{[V^k]}$ 
lie on an extremal face of $\cone(X)$, or $\pi^*H$ is
negative on an extremal ray, whose locus has to be contained in
the indeterminacy locus of $\pi$ and therefore has codimension
greater than one in $X$. 
\end{proof}


\section{Examples}

In this section we show the effectiveness of all cases listed in
theorem \ref{main}. If $X$ has only fiber type contractions,
examples are given by the products such as $\proj^{i_X-1} \times Y$, where
$Y$ is a suitable fourfold of pseudoindex $i_X$ (for the ones with $i_X=2$ see \cite{Wifour}). 
The remaining cases are listed below:\par

\medskip
\begin{center}
\begin{tabular}{|c|c||c|c|c|c||c|}\hline
\quad $\dim X$ \quad &\quad $\rho_X$ \quad  & \quad $R_1$ \quad & \quad $R_2$ \quad & \quad
$R_3$ \quad & \quad $R_4$ \quad & \quad \quad \\
\hline\hline
 $5$ & $2$  & $F$ & $D_0$  & & & a \\
   &   & $F$ & $D_1$ & & & b \\
   &  & $F$ & $D_2$  & & & c \\
   &   & $F$ & $S$  & & & d \\
   &   & $D_2$ & $D_2$  & & & e \\
&	  & $D_2$ & $S$  & & & f \\
	  \hline
& $3$  & $F$ & $F$ & $S$ & & g \\
 &     & $F$ & $F$ & $D_1$ & & h \\
 &     & $F$ & $F$ & $D_2$ & & i \\
 &     & $F$ & $D_2$ & $D_2$ & & j \\
	  \hline
& $4$  & $F$ & $F$ & $F$ & $D_2$ & k \\
\hline\hline
$6$ & $2$  & $F$ & $D_1$  & & & l \\
 &   & $F$ & $D_2$  & & & m \\
 &   & $F$ & $S$  & & & n \\
\hline\hline
$7$ & $2$  & $F$ & $D_2$  & & & o \\
\hline
\end{tabular}
\end{center}

\bigskip
\begin{itemize}
\item[\bf{a.}] \qquad $X = \proj_{\proj^4}(\Ol \oplus \Ol(1)) = Bl_p\proj^5$.\par

\medskip
\item[\bf{b1.}]\qquad $X = \proj_{\pt}(\Ol \oplus \Ol \oplus\Ol(2))$.

\item[\bf{b2.}]\qquad $X = \proj_{\quadr^3}(\Ol \oplus \Ol \oplus \Ol(1)).$\par

\medskip
\item[\bf{c.}] \qquad $X = \proj_{\pd}(\Ol \oplus \Ol \oplus \Ol \oplus \Ol(1)) 
= Bl_{\pd} \proj^5$.\par

\medskip
\item[\bf{d.}] \qquad $X = \proj_{\pt}(\Ol \oplus \Ol(1) \oplus \Ol(1))$.\par

\medskip
\item[\bf{e1.}] \qquad Let $X=Bl_{S_3}\proj^5$, where $S_3$ is a cubic scroll contained in a hyperplane 
$H \subset \proj^5$; denote by  $\sigma$ the blow-up and by $E$ the exceptional divisor.\\ 
Let $\sigma^*\Ol_{\proj^5}(1)$ be the pull-back to $X$ of the hyperplane bundle of $\proj^5$, and 
let $\acca = \sigma^*\Ol(1) -E$ be
the strict transform of $H$; the linear system
$$|L|=\sigma^*|\Ol(2) - S_3|=|2\sigma^*\Ol(1) -E|$$ 
has empty base locus on $X$ 
and the associated map $\f_{|L|}$ gives $\acca$ a structure of $\pdue$-bundle over $\pdue$.\\
Moreover $\acca_{|\acca}=(L-\sigma^*\Ol(1))_{|\acca}$, so that the restriction of $\acca$ to each
fiber of $\f_{|L|}$ is $\Ol_{\pdue}(-1)$; we can therefore apply the Nakano contractibility 
criterion \cite{Na}, which yields the existence of a manifold $M \supset \pdue$ such that 
$X \simeq Bl_{\pdue}(M)$ and $\acca$ is the exceptional divisor of this blow-up.\\ 
Moreover, if we denote by $\psi$ the rational map associated to the li\-near system 
$|\Ol(2) - S_3|$ on $\proj^5$ we have that the following diagram commutes:

$$\xymatrix@=30pt{ & X \ar[ld]_\sigma \ar[rd]^{\f_{|L|}}& \\
\proj^5   \ar@{-->}[rr]^\psi&  & M}$$

One can also prove (see \cite{Fu1}) that $M$ is isomorphic to the hyperplane section 
associated with the Pl\"ucker embedding of the Grassmannian 
$G(1,4) \subset \proj^9$, so (see [7.1] in \cite{Fu1}) $M$ is a del Pezzo variety.  

Note that since $\rho_X = 2$ and $X$ has two smooth blow-downs, 
$-K_X$ is positive on the entire cone $\cone(X)$, so $X$ is a Fano variety. 
Moreover, we can write 
   $$-K_X = 6 \sigma^*\Ol(1) -2E = 2(3 \sigma^*\Ol(1)- E),$$
so $X$ has index $2$.\par

\medskip
\item[\bf{e2.}] \qquad Let $X=Bl_{V}\proj^5$, where $V$ is a Veronese surface.
Denote by $\sigma$ the blow-up and by $E$ the exceptional divisor.\\
Consider on $\proj^5$ the linear system $|\Ol_{\proj^5}(2)-V|$  of the quadrics containing $V$, and 
denote by $\xymatrix@1{\proj^5~ \ar@{-->}[r]^F& ~\proj^5}$ the associated rational map; 
call $V'$ the exceptional locus of
$F^{-1}$ and let $\sigma': X' \to \proj^5$ be the blow-up of $\proj^5$ along $V'$.

$$\xymatrix{ X \ar[d]_\sigma	&	\cong   & 	X' \ar[d]^{\sigma'}\\
\proj^5	\ar@{-->}[rr]^F &	  &	\proj^5}$$

One can prove (see [2.0.2] in \cite{ESB}) that $X' \simeq X$, that the exceptional 
divisors of the two blow-ups satisfy the relations
\begin{itemize}
\item[] $E = 2\sigma^*\Ol(1) - \sigma'^*\Ol(1)$
\item[] $E' = 2\sigma'^*\Ol(1) - \sigma^*\Ol(1),$
\end{itemize}
and that the map $F$ is an involution (Theorem 2.6 in \cite{ESB}).\par
As in the previous example, since $\rho_X = 2$ and $X$ has two smooth blow-downs, 
$-K_X$ is positive on the entire cone $\cone(X)$, so $X$ is a Fano variety, and
from the canonical bundle formula 
   $$-K_X = 6 \sigma^*\Ol_{\proj^5}(1) -2E = 2(3\sigma^*\Ol_{\proj^5}(1) - E),$$
we have that $X$ has index $2$.\par

\medskip
\item[\bf{f.}] \qquad Let $X = Bl_{\quadr^2}\proj^5$, where $\quadr^2 \subset \proj^5$
is a smooth two-dimensional quadric, denote by $\sigma$ the blow-up and by $E$ the
exceptional divisor; then
   $$-K_X = 6 \sigma^*\Ol_{\proj^5}(1) -2E.$$
For every curve $C \subset X$ which is not contained in $E$, we have that $\sigma(C)$
is a curve in $\proj^5$ of a certain degree $d$, and the sum of the multiplicities
of the points of intersection of $\sigma(C)$ and $\quadr^2$ is $\le 2d$. This implies that 
   $$-K_X \cdot C \ge 6d - 4d \ge 2d.$$
The exceptional divisor $E$ can be written as 
   $$E = \proj(\Norm^*_{\quadr^2|\proj^5}) = \proj(\Ol(-2) \oplus \Ol(-1) \oplus \Ol(-1)),$$ 
and $E_{|E} = - \xi_{\Norm^*}$.
If we denote by $Q$ the section of $\sigma:E \to \quadr^2$ which corresponds to $\Ol(-2)$, 
then $\cone (E) = \threespan{[l_1]}{[l_2]}{[l_3]}$, where $l_1$ and $l_2$ correspond to the 
two rulings of $Q$ and $l_3$ is a line in a fiber of $\sigma$.\\
If we write $K_E$ as $-3\xi -6\sigma^*\Ol_{\quadr^2}(1)$, the adjunction formula yields
   $$-K_{X_{|E}} = -K_E + E_{|E} = 2\xi + 6\sigma^*\Ol(1),$$
so $-K_X \cdot l_i = 2$ for every $i$, hence $X$ is a Fano variety of pseudoindex $2$.\par

\medskip
The line bundle $2\sigma^*\Ol(1) - E$ is nef on $X$, and it vanishes on the strict transform
of the $\pt \subset \proj^5$ which contains $\quadr^2$; hence it is the supporting divisor of
the small contraction of $\pt$ to a point.\par

\medskip
\item[\bf{g.}] \qquad Let $\E = \Ol_{\proj^3} \oplus \Ol_{\proj^3}(1)$
and $Y = Z = \proj_{\proj^3}(\E) = Bl_p \proj^4$;
consider the fiber product $X= Y \times_\pt Z$

$$\xymatrix@=30pt{
X \ar[d]_{p_Z} \ar[r]^{p_Y} & Y \ar[d]^\f \ar[r]^\sigma & \proj^4\\
Z \ar[r] & \proj^3 &
}$$

\noindent Call
\begin{itemize}
\item[] $\xi_Y \in \pic(Y)$ the tautological bundle of $\E$,
\item[] $H=\f^*(\Ol_{\proj^3}(1))$,
\item[] $E$ the exceptional divisor of $\sigma$,
\item[] $L=\sigma^*(\Ol_{\proj^4}(1))$.
\end{itemize}

\noindent Then
$$\pic(Y)=\langle \xi_Y,H \rangle = \langle E,L \rangle$$
and the canonical bundle of $Y$ can be written as
$$K_Y  = - 2\xi_Y + \f^*(K_{\proj^3} + \det \E) = - 2\xi_Y -3H$$
or
$$K_Y = \sigma^*K_{\proj^4} + 3E = -5L + 3E.$$
Note also that $\xi_Y= L$, and since for every fiber $f$ of $\f$
we have $H \cdot f = 0$, $E \cdot f = L \cdot f = 1$ and for
every line $l \subset E$ which is not contracted by $\sigma$ we
have $H \cdot l = 1$, we can also write $H = L-E$.\par

\medskip
Now consider on $Y$ the rank 2 vector bundle $\F = \Ol_Y \oplus
H$: then $X$ can be seen as $\proj_Y(\F)$. 
Call $\widetilde Y$ the section of $p_Y$ which corresponds to the surjection $\F \to \Ol_Y \to 0$. 
If we denote by $\xi$ the tautological bundle of $\F$ we can write
\begin{eqnarray*}
-K_X  & = & 2\xi - p_Y^*(K_Y + \det \F)\\
       & = & 2(\xi + p_Y^*(\xi_Y + H)),
\end{eqnarray*}
or, in terms of $L$ and $E$,
$$-K_X = 2(\xi + p_Y^*(2L-E)).$$

\medskip
First of all we show that $X$ is a Fano variety with $r_X = 2$.\\
The line bundle $2L - E$ is ample on $Y$, so $p_Y^*(2L-E)$ 
is nef on $X$; since also $\xi$ is nef,
we have to show that on every curve in $X$ at least one of these bundles is nonzero.\\
But the line bundle $p_Y^*(2L-E)$ vanishes only on the fibers of $p_Y$, where
we know that $\xi$ is positive, so $X$ is a Fano variety of index $2$.\par

\medskip
$X$ has two fiber type contractions, associated with the nef bundles
$p_Y^*(2L-E)$ and $p_Y^*H + \xi$.\\
The line bundle $\xi + p_Y^*\xi_Y$ supports a small contraction which
contracts a $\proj^3 \subset \widetilde Y$ to a point. 
In this case there exist two divisors on $X$ which are negative on the exceptional locus of this contraction: 
one is $\widetilde Y$ and the other is $p_Y^*E$: in fact, since
$E = L-H$, we have $p_Y^*E \cdot l = -1$ for every line $l \subset \proj^3$. \par

\medskip
\item[\bf{h1.}] $\qquad X = \proj_{\proj^3}(\Ol \oplus \Ol(2)) \times \pu$.
\item[\bf{h2.}] $\qquad X =\proj_{\proj^3}(\Ol \oplus \Ol(1)) \times \pu$.
\item[\bf{h3.}] $\qquad X = \proj_{\quadr^3}(\Ol \oplus \Ol(1)) \times \pu$.\par

\medskip
\item[\bf{i.}] \qquad $X = Bl_l \proj^4 \times \pu$.\par

\medskip
\item[\bf{j1.}] \qquad Here we construct an example of a fivefold where the two 
divisorial contractions have the same exceptional locus.

Let $\F = \Ol_{\pdue \times \pdue} \oplus \Ol_{\pdue \times \pdue}(1,1)$ and 
$X = \proj_{\pdue \times \pdue}(\F)$,
let $\pi$ be the projection map, $\xi$ the tautological bundle on $X$ and 
$E$ the section of $\pi$ which corresponds to the surjection
$\F \to \Ol_{\pdue \times \pdue} \to 0$. Then
$$-K_X = 2\xi - \pi^*(K_{\pdue \times \pdue} + \det \F) = 2(\xi + \pi^*\Ol(1,1));$$
since $\pi^*\Ol(1,1)$ vanishes only on the fibers $f$ of $\pi$, while $\xi \cdot f = 1$, it follows that 
$X$ is a Fano variety and $i_X = 2$.\par

\smallskip
Obviously $X$ admits a fiber type contraction, which is given by
its structure of $\pu$-bundle and is supported by $\pi^*\Ol(1,1)$.\\
The nef bundles $\xi + \pi^*\Ol(1,0)$ and $\xi + \pi^*\Ol(0,1)$ vanish each on one ``ruling" of $E$,
so they support two different contractions of $E$ to $\pd$, 
which are in fact smooth blow-downs.\\
Finally, the line bundles
$\pi^*\Ol(1,0)$ and $\pi^*\Ol(0,1)$ support the contractions of
the two faces of $\cone(X)$ which contain the fiber type ray, as
shown in the diagram.

$$
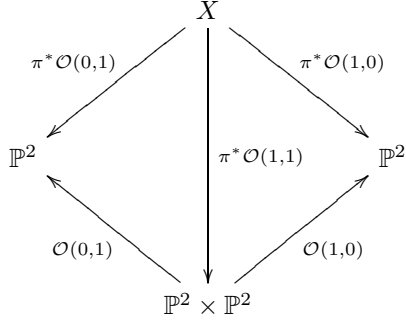
\xymatrix@=40pt{ & X\ar[dd]^{\pi^*\Ol(1,1)} \ar[ld]_{\pi^*\Ol(0,1)} \ar[rd]^{\pi^*\Ol(1,0)}& \\
\pdue & & \pdue\\
& \pdue \times \pdue  \ar[lu]^{\Ol(0,1)}\ar[ru]_{\Ol(1,0)} &}$$

\medskip
\item[\bf{j2.}]\qquad An example of a fivefold with two divisorial contractions with disjoint 
exceptional loci is $X=Bl_{\Pi_1 \sqcup \Pi_2} \proj^5$,
the blow-up of $\proj^5$ along two disjoint planes $\Pi_1$, $\Pi_2$.\\
Let $\sigma:X \to \proj^5$ be the blow-up and denote by $E_1$ and $E_2$
the exceptional divisors; by the canonical bundle formula we have
$$-K_X = 6\sigma^*\Ol_{\proj^5}(1)-2E_1-2E_2=2(3\sigma^*\Ol_{\proj^5}(1)-E_1-E_2).$$
We want to prove that $H:=3\sigma^*\Ol_{\proj^5}(1)-E_1-E_2$ is ample on $X$:
if $C$ is a curve not contained in $E_1 \cup E_2$, then $\sigma(C)$ is a curve of degree $d$ in $\proj^5$
which intersects $\Pi_1$ and $\Pi_2$ in a number of points which has to be
less or equal than $d$ (counted with multiplicity).
So
$$H \cdot C \ge d \ge 1.$$
As for curves contained in an exceptional divisor $E_i$, we know that
$$E_i=\proj_\pd(\Norm^*_{\Pi_i|\proj^5})=\proj(\Ol(-1)^{\oplus 3})\simeq \pdue \times \pdue$$
so
$$-K_{E_i}=3 \xi_i + 6\pi^*\Ol_\pd(1),$$
where $\xi_i$ is the tautological bundle of $\Norm^*_{\Pi_i|\proj^5}$ and $\pi$ is the projection $E_i \to \pd$.
By the adjunction formula, recalling that ${-E_i}_{|E_i}=\xi_i$ we have that
$$(-K_X)_{|E_i}= -K_{E_i} + {E_i}_{|E_i}=3 \xi_i + 6\pi^*\Ol(1) - \xi_i = 2 \xi_i + 6\pi^*\Ol(1),$$
so $H_{|E_i} = \xi_i + 3\pi^*\Ol(1) \simeq \Ol_{\pd \times \pd}(1,-1) \otimes 
\Ol_{\pd \times \pd}(0,1)^{\otimes 3} \simeq \Ol_{\pd \times \pd}(1,2)$ which is ample, 
hence we have proved that $X$ is a Fano variety of index $2$.\\

Now consider in $\proj^5$ the lines which intersect both $\Pi_1$ and $\Pi_2$, and call $V$ the family
of deformations of their strict transforms on $X$.
Then $V$ is covering (since lines meeting $\Pi_1$ and $\Pi_2$ cover the whole $\proj^5$) 
and $-K_X \cdot V = 2$, so $V$ is unsplit. Hence $V$ is extremal and is associated to a fiber type contraction
$\f:X \to Y$, which can be easily proved to be a $\pu$-bundle over a smooth fourfold.\\
The fibers of $\f$ are the strict transforms of the lines meeting $\Pi_1$ and $\Pi_2$,
hence, for any fiber $f$ we have $E_i \cdot f=1$; being $E_1 \cap E_2=\emptyset$ it follows
that $E_1$ and $E_2$ are sections of $\f$.\\
It's now easy to prove that $Y \simeq \pd \times \pd$ and $X=\proj_Y(\Ol(1,2)\oplus \Ol(2,1))$.\par

\medskip
\item[\bf{j3.}]\qquad In this example the exceptional loci of the two divisorial 
contractions of $X$ are different 
but have nonempty intersection.\\
Let  $Y = Bl_l \proj^4 = \proj_{\pd}(\Ol \oplus \Ol \oplus \Ol(1))$.

$$\xymatrix@=30pt{Y  \ar[r]^\sigma\ar[d]_\f & \proj^4 \\
\proj^2 &}$$

Call $H=\f^*\Ol_\pd(1)$ and let
$$X=\proj_Y(H \oplus \Ol_Y(1)).$$

Using the same notation as in example {\bf g.} we can write
$$-K_X = 2(\xi + \pi^*(\xi_Y + H)),$$
so again we have that $X$ is a Fano variety of (pseudo)index 2.

$X$ admits a fiber type contraction on $Y$, which is supported by
$\pi^*(L + H)$, and two divisorial contractions: the first one
contracts the special section on $X$ to $\pd$ and is supported by $\xi + \pi^*H$, while the second one
contracts the $\pu$-bundle over the exceptional divisor in $Y$ to a two-dimensional quadric,
and is supported by $\xi + \pi^*L$.\par

\medskip
\item[\bf{k.}] \qquad $X = \pu \times \pu \times Bl_p\pt$.\par

\medskip
\item[\bf{l.}] \qquad $X = \proj_{\proj^4}(\Ol \oplus \Ol \oplus \Ol(1)) = Bl_l\proj^6$.\par

\medskip
\item[\bf{m.}] \qquad $X = \proj_{\proj^3}(\Ol \oplus \Ol \oplus\Ol \oplus \Ol(1)) 
= Bl_{\pdue}\proj^6$.\par

\medskip
\item[\bf{n.}] \qquad $X = \proj_{\proj^4}(\Ol \oplus \Ol(1) \oplus \Ol(1))$.\par

\medskip
\item[\bf{o.}] \qquad $X = \proj_{\proj^4}(\Ol \oplus\Ol \oplus\Ol \oplus \Ol(1)) 
= Bl_{\pdue}\proj^7$.\par

\end{itemize}


\section{Fano fivefolds with a covering quasi-unsplit family} \label{with}

In this section we start the proof of theorem \ref{main}.
We know as a general fact (see subsection \ref{fano}) that on $X$ there exist 
covering locally unsplit families of rational curves of degree $\le \dim X + 1$; 
since we are assuming $\rho_X \ge 2$ these families have degree $\le 5$, otherwise 
proposition \ref{iowifam} and corollary \ref{rho1} would imply $\rho_X = 1$.\\
We start considering the case when one of these families is quasi-unsplit.
Note that this is the case if on $X$ there exists a fiber type ray $R$: in fact, in this case, 
through every $x \in X$ there exists a rational curve which is contracted by $\f_R$ and 
has degree $\le \dim X + 1$; among the families of deformations of these curves we can choose 
a covering one with minimal degree, which is quasi-unsplit since $R$ is extremal.\par

\medskip
\begin{lemma} \label{3uns}
Let $V$ be a covering quasi-unsplit locally unsplit family of rational curves and
$R_1$ an extremal ray of $\cone(X)$ independent from $[V]$; 
assume that the contraction $\f_{R_1}$ has a three-dimensional fiber $F$.
Then there exists a covering unsplit family which is numerically proportional to $V$.
\end{lemma}

\begin{proof}
If $V$ is not unsplit then $\deg V \ge 4$, so $V_x$ cannot be unsplit
for any $x \in F$, otherwise we would have $\dim(\loc(V_x) \cap F) \ge 1$ against the fact
that $\cone(\loc(V_x))=\onespan{[V]}$ and $\cone(F)=\onespan{R_1}$. \\
Therefore through every point of $F \cap \loc(V)$ there passes a reducible cycle
in $\V$. Moreover, $\loc(\V)$ is closed since $\V$ is proper, and since $V$ 
is covering $\loc(\V)=X$, hence through any point of $F$ there is a reducible cycle in 
$\V$.\\
It follows that $F$ is contained in the
locus of the family of deformations of one of the components of these cycles.
Note that, since $\deg V \le 5$ and $V$ is quasi-unsplit, such a family is 
unsplit and numerically proportional to $V$.\\
We denote this family by $\alpha V$, and applying lemma \ref{locy}
(a) we obtain that $\dim \loc(\alpha V)_F \ge 5$, so $\alpha V$ is covering.  
\end{proof}

\smallskip
\begin{lemma} \label{triangle}
Let $V$ be a covering unsplit family of rational curves and
$R_1$ an extremal ray of $\cone(X)$ independent from $[V]$. 
Assume that the contraction $\f_{R_1}$ has a three-dimensional fiber $F$ and let 
$D$ be an irreducible component of $\loc(V)_F$ (note that, by lemma \ref{locy},
$\dim D \ge 4$). Then
   \begin{itemize}
   \item[(a)] if either $D = X$ or $D \cdot V > 0$ then $\cone(X) = \twospan{[V]}{R_1}$;
   \item[(b)] if $R_2$ is a birational ray different from $R_1$ then $D \cdot R^2 = 0$;
   \item[(c)] if $R_2$ is a divisorial ray and $E_2$ is its exceptional locus, then\\
   $E_2 \cdot V = E_2 \cdot R^1 = 0$;
   \item[(d)] if $R_2$ is a fiber type ray then $D \cdot R^2 > 0$.
   \end{itemize}
\end{lemma}

\begin{proof} 
The proof of (a) is an easy consequence of
corollary \ref{ray}: in fact, we know that $\cone(D)=\twospan{[V]}{R_1}$,
and if $D$ is a divisor and $D \cdot V > 0$ we can write $X = \cloc_2(V)_F$ (since $V$
is covering), while if $D=X$ the proof is trivial.\\
To prove (b), we observe that the nontrivial fibers of $\f_{R_2}$ 
have dimension $\ge 2$, so if $D
\cdot R^2 \neq 0$ then $D$ contains a curve whose numerical
class is in $R_2$, a contradiction.\\
In case (c), if either $E_2 \cdot V > 0$ or $E_2 \cdot R^1 \neq 0$ then
$E_2 \cap D \neq \emptyset$.
Take a point $x \in E_2 \cap D$ and a curve $C$ in $R^2$ passing
through $x$: by (c) we know that $D \cdot C = 0$, so $C \subseteq D$,
a contradiction.\\
Finally, to prove (d) let $x$ be a point in $D$ and $C$ a curve in $R^2$ through $x$.
Since $C$ cannot be contained in $D$ we must have $D \cdot R^2 > 0$. 
\end{proof}

\smallskip
\begin{remark}\label{triangle'}
Note that (b), (c) and (d) still hold if we replace $R^2$ with any
noncovering unsplit family which is independent from $V$ and $R^1$.
\end{remark}

\smallskip
\begin{lemma} \label{negray}
Let $R$ be a divisorial ray on $X$, let $E$ be its
exceptional locus and consider the intersection number $E
\cdot R^i$ with all the divisorial extremal rays of $\cone(X)$
different from $R$.
Then $E \cdot R^i < 0$ for at most one index $i$; moreover
in this case $E \cdot R^j =0$ for every $j \neq i$ and $\cone(E) = \twospan{R}{R_i}$.
\end{lemma}

\begin{proof}
Assume that there exists an index $i$ such that $E \cdot R^i <0$;
then we have $E = \loc(R)_{\loc(R^i_x)}$, so $\cone(E) = \twospan{R}{R_i}$ by corollary \ref{ray}. 
In particular, $E$ cannot contain curves whose class is in $R_j$ for $j \neq i$.
\end{proof}

\bigskip
We can now start the proof of theorem \ref{main}.

\subsection*{\boldmath{$\rho_X = 2$}\unboldmath}

\quad We have to prove that at least one of the extremal contractions on $X$ is
of fiber type. Assume that this is not the case; in particular
$[V]$ is not extremal, so by lemma \ref{codim}
there exists a small extremal ray $R_1$.\\
Denote by $R_2$ the other extremal ray of $\cone(X)$; by lemma \ref{3uns}
we can assume that $V$ is unsplit, and by lemma
\ref{triangle} either $\cone(X) = \twospan{[V]}{R_2}$ and $[V]$ is
extremal or there exists an effective divisor
$D$ such that $D \cdot V = D \cdot R^2 = 0$, implying that $D$ is
numerically trivial on $\cone(X)$; in both cases we reach a contradiction.


\subsection*{\boldmath{$\rho_X = 3$}\unboldmath}

\quad We divide this part of the proof into three cases.\par
\medskip
{\bf Case 1.} \quad All rays of $\cone(X)$ are of fiber type.\par
\smallskip
If two rays, say $R_1$ and $R_2$, do not lie on the same
extremal face of $\cone(X)$, we can consider the rationally
connected fibration $\pi:\xymatrix@1{X~ \ar@{-->}[r]& ~Z}$ associated to $R^1$ and $R^2$. 
Since $\rho_X = 3$ we have $\dim Z > 0$, so by
lemma \ref{codim} $X$ must have a small elementary
contraction, a contradiction. The only possibility to exclude this
situation is that $\cone(X)$ has exactly three rays.\par

\medskip
{\bf Case 2.} \quad In $\cone(X)$ there exists a small extremal ray.\par
\smallskip
In this case we prove that $\cone(X)=\threespan{R_1}{R_2}{R_3}$, where $R_1$ is small
and both $R_2$ and $R_3$ are of fiber type.\\
Denote the small ray by $R_1$, and denote by $F_1$ an irreducible component of a fiber
of $\f_{R_1}$. Note that by lemma \ref{3uns} we can assume that $V$ is unsplit.

First of all we prove that $X$ has at least one fiber type contraction:
suppose that this is not the case, let $D_1= \loc(V)_{F_1}$ and
apply lemma \ref{triangle}. Since $\rho_X=3$ we cannot be in case
(a), and so $D_1$ is a divisor such that
$D_1 \cdot R^i = 0$ for every $i \neq 1$; as a consequence
$\cone(X) = \threespan{R_1}{R_2}{R_3}$.\\
If $R_2$ is a small ray, we can repeat the same argument
with the divisor $D_2 = \loc(V)_{F_2}$, and we obtain that $D_2$
vanishes on the face $\twospan{R_1}{R_3}$; since $D_1$
vanishes on the face $\twospan{R_2}{R_3}$ and $D_1 \cdot V = D_2 \cdot V
=0$, it must be
$[R^3] = [V]$, against the assumption that $X$ has no fiber type contractions.\\
So both $R_2$ and $R_3$ are divisorial. By lemma \ref{triangle}
(c), if we denote by $E_i$ the exceptional locus of $R_i$ we have
$E_i \cdot V = E_i \cdot R^1 = 0$, and we know that $E_i \cdot R^i
< 0$, which implies $E_2 \cdot R^3 >0$ and $E_3 \cdot R^2 > 0$;
in particular this yields that the intersection
numbers of $E_2$ and $E_3$ with every curve in $X$ have opposite signs.
The existence of curves which intersect $E_2 \cup E_3$ without
being contained in it gives rise to a contradiction.\\
We have thus proved  that $X$ has at least one fiber type contraction, 
associated to a ray $R_2$. \par

\smallskip
Suppose by contradiction that every other ray $R_i$ of $\cone(X)$ is birational.
By lemma \ref{triangle} (b), for the divisor $D_1^2:=\loc(R^2)_{F_1}$ 
we have $D_1^2 \cdot R^i = 0$; moreover lemma \ref{triangle} (a) implies that
$D_1^2 \cdot R^2 = 0$, so $\cone(X) = \threespan{R_1}{R_2}{R_3}$.\\
The ray $R_3$ cannot be divisorial, otherwise we
would have by lemma \ref{triangle} (c) that $E_3 \cdot R^1 = E_3
\cdot R^2 =0$ while $E_3 \cdot R^3 < 0$, against the
effectiveness of $E_3$, so $R_3$ must be small.\\
Let $F_3$ be an irreducible component of a fiber of $\f_{R_3}$ and let
$D_3^2:=\loc(R^2)_{F_3}$; by lemma \ref{triangle} we have
$D_3^2 \cdot R^1 = D_3^2 \cdot R^2 =0$.\\
Consider a minimal horizontal dominating family $V'$ for the fiber
type contraction $\f_{R_2}$; from the results in Section 8, Case 1, \cite{ACO}
we know that $V'$ is unsplit.\\
The family is independent either from $R_1$ and $R_2$ 
or from $R_2$ and $R_3$; assume without loss of generality that we are in the
first case.\\
If $V'$ is covering we have $X = \loc(V',R^2)_{F_1} = \loc(R^2,V')_{F_1}$, so
$R^3=V'$ and $R_3$ is of fiber type, a contradiction.\\
If else $V'$ is noncovering, then by remark \ref{triangle'} we have
$D_1^2 \cdot V'= D_3^2 \cdot V'=0$, so that $[R^2]=[\lambda V']$,
again a contradiction.\par
\smallskip
We have thus proved that $X$ admits a small ray $R_1$ and at least two fiber type
rays $R_2$, $R_3$; by lemma \ref{3uns} we know that the families
$R^2$, $R^3$ are unsplit, so by lemma \ref{locy}
we have that $X = \loc(R^3,R^2)_{F_1} = \loc(R^2,R^3)_{F_1}$ and lemma \ref{ray}
implies that $\cone(X) = \threespan{R_1}{R_2}{R_3}$.\par

{\bf Case 3.} \quad In $\cone(X)$ there is at least a birational ray, but no small rays.\par
\smallskip
In this case we prove that $\cone(X)=\threespan{R_1}{R_2}{R_3}$, where at least
one $R_i$ is of fiber type, and that the possible cases are the ones listed in theorem \ref{main}.\\
Since $X$ has no small contractions we know by lemma \ref{codim}
that $[V]$ lies on an extremal face of $\cone(X)$.\par

\smallskip
Suppose that there exists a ray $R_1$ which does not lie in a face with
$[V]$, and denote by $E_1$ its exceptional locus.\\
If either $R_1$ is divisorial and $E_1 \cdot V > 0$ or
$R_1$ is of fiber type then the associated family $R^1$
is horizontal and dominating with respect to the rc$V$-fibration.
Hence we can apply lemma \ref{codim} to $V$ and $R^1$ and conclude
that $[V]$ and $[R^1]$ are on the same extremal face, a contradiction.\\
So we can assume that $R_1$ is divisorial and $E_1 \cdot V = 0$.
Then $E_1$ must be negative on another ray $R_2$ which lies in a face with $[V]$: in fact, 
$E_1$ cannot vanish on a face containing $[V]$, otherwise it would be $\le 0$
on the entire cone; clearly $R_2$ has to be divisorial. Then we can
conclude from lemma \ref{negray} that in $\cone(X)_{E_1<0}$ there are 
two divisorial rays, in $\cone(X)_{E_1>0}$ there are only fiber type rays
and $\cone(E_1)=\twospan{R_1}{R_2}$.\\
Let $R_3$ be one of the fiber type rays; we can write $X=
\loc(R^3)_{E_1}$, and we have by
remark \ref{face} that $\cone(X) = \threespan{R_1}{R_2}{R_3}$. \par

\smallskip
In the case when every extremal ray lies on
a face with $V$ we have trivially that $\cone(X) = \threespan{[V]}{R_1}{R_2}$.\par

\smallskip
If $X$ has two fiber type rays $R_1$, $R_2$ and one divisorial ray $R_3$, then
$\f_{R_3}$ cannot have a four-dimensional fiber $F_3$, otherwise we would have
$X= \loc(R^1)_{F_3}$ and $\rho_X=2$, by remark \ref{face}.\\
Finally, in the case when $X$ has one fiber type ray $R_1$ and two divisorial rays
$R_2$ and $R_3$, we claim that both $R_2$ and $R_3$ have
two-dimensional fibers: in fact, if $R_2$ has a fiber $F_2$ of
dimension three, by lemma \ref{3uns} and lemma \ref{triangle} (c) we have that
$E_3 \cdot V = E_3 \cdot R^2 = 0$, a contradiction. 

\subsection*{\boldmath{$\rho_X = 4$}\unboldmath}
\quad In this case (see Section 8, Case 2 in \cite{ACO}) $X$ is rationally connected 
with respect to four independent unsplit families $V$,$V'$, $V''$ and $V'''$, such
that each one is horizontal with respect to the fibration
associated to the previous ones.\\
By remark \ref{rcvdim}, for three
among these families, say $V$, $V'$ and $V''$, the pointed locus has dimension 1, 
so these families are covering.\\
Moreover, if there exists a small ray $R$ we can
choose two covering families, say $V$ and $V'$, such that $[V],
[V']$ and $[R]$ are numerically independent; then if $F$ is a
fiber of $\f_R$ we can write $X = \loc(V, V')_F$,
implying that $\rho_X = 3$, a contradiction. So two cases are
possible: either all rays are of fiber type or there
exists a divisorial ray. \par

\smallskip
Suppose that all the rays of $\cone(X)$ are of fiber type. If 
there exist two rays $R_1$, $R_2$ which do
not lie on the same extremal face of $\cone(X)$, we can consider the rationally
connected fibration $\pi:\xymatrix@1{X~ \ar@{-->}[r]& ~Z}$ associated to $R^1$ and $R^2$. 
Since $\rho_X = 4$ we have $\dim Z > 0$, so by
lemma \ref{codim} $X$ must have a small elementary
contraction, a contradiction.\\
So every pair of extremal rays lies on an extremal two-dimensional face of
$\cone(X)$; it is easy to verify that in this case $\cone(X)$ has
exactly four rays. \par

\smallskip
Suppose now that there exists a divisorial ray $R$.\\
Since $X$ has no small contractions and $\rho_X=4$, $V, V'$
and $V''$ lie on the same extremal face $\sigma$ of $\cone(X)$ by
lemma \ref{codim}, and, applying again lemma \ref{codim} to every
pair of families chosen among $V, V'$ and $V''$, we get that
$\sigma = \threespan{[V]}{[V']}{[V'']}$.\\
Let $F$ be a fiber of $\f_{R}$, which has dimension greater than two by proposition
\ref{fiberlocus}. Since $R$ is not in
$\sigma$ we have $\dim \loc(V,V',V'')_F \ge \dim F + 3$ by lemma
\ref{locy}, so $\dim F=2$, $X = \loc(V,V',V'')_F$ and every curve in $X$ can
be written with positive coefficients with respect to $R$ and $V$;
but $V$, $V'$ and $V''$ play a symmetric role, so we can conclude
that $\cone(X) = \langle [V], [V'], [V''], [R] \rangle$.


\section{Fano fivefolds without a covering quasi-unsplit family} \label{without}

In this section we conclude the proof of theorem \ref{main} considering Fano
fivefolds $X$ which do not have any covering quasi-unsplit locally
unsplit family; more precisely we prove the following \par

\smallskip
\noindent{\bf Theorem \ref{shelob}.}
\emph{Let $X$ be a Fano fivefold of pseudoindex two which does not have a 
covering quasi-unsplit locally unsplit family of rational curves; then $\rho_X = 2$ and $X$ 
is  the blow-up of $\proj^5$ along a two-dimensional smooth quadric (a section of $\Ol(2)$
in a linear $\proj^3 \subset \proj^5$), or along a cubic scroll ($\proj_\pu(\Ol(1) \oplus \Ol(2))$ 
embedded in an hyperplane of $\proj^5$ by the tautological bundle), or along a Veronese surface.}\par

\smallskip
\begin{proof}
Let $V$ be a locally unsplit dominating family on $X$ and let $\V$ be the associated Chow family.
Since $V$ is not quasi-unsplit then $[V]$ cannot be extremal; in particular it follows that $\rho_X \ge 2$.
Moreover, since $V$ is locally unsplit but not unsplit we have 
$$4 = 2i_X \le \deg V \le \dim X +1 = 6;$$ 
moreover, if $\deg V = 6$ then we would have $X = \loc(V_x)$ for a general $x \in X$ and 
$\rho_X = 1$ by corollary \ref{rho1}, a contradiction,
therefore we can assume that $4 \le \deg V \le 5$.
As a consequence we have that every reducible cycle 
in $\V$ splits into exactly two irreducible components.\par

\medskip
Consider the pairs $(W^i,\W^i)$ of families such that
there is a cycle in $\V$ whose irreducible components belong respectively to $W^i$ and $\W^i$,
and let $\B$ be the set of these pairs. By this definition we clearly have
$[W^i]+[ \W^i] = [V]$, and since the anticanonical degree of these families is bounded they are only
a finite number.\\
We begin establishing some properties of these pairs.\par

\smallskip
\begin{lemma}\label{dimwi}
If $(W^i,\W^i) \in \B$ then the families $W^i$, $\W^i$ are unsplit, and moreover\linebreak 
$(\dim \loc(W^i),\dim \loc(W^i_x))$ is either $(4,2)$, $(4,3)$, $(4,4)$ or $(3,3)$.
\end{lemma}

\begin{proof} 
The families are unsplit since 
\begin{equation*}
4 = 2i_X \le \deg W^i + \deg \W^i = \deg V \le 5,
\end{equation*}
and therefore they are noncovering, so the second assertion follows from proposition
\ref{iowifam}. 
\end{proof}

\smallskip
\begin{lemma}\label{dimdi}  
Let $(W^i,\W^i) \in \B$ such that $[W^i] \neq [\alpha V]$, and let $D_i$ and $\Di_i$ be meeting 
components of $\loc(W^i)$ and $\loc(\W^i)$. Up to exchange $D_i$ and $\Di_i$, we have that
$(\dim D_i, \dim \Di_i)$ is either $(3,4)$ or $(4,4)$.
\end{lemma}

\begin{proof} 
By lemma \ref{dimwi} we have $\dim D_i$, $\dim \Di_{i} \ge 3$, and
equality holds if and only if $D_i = \loc(W^i_x)$ for some $x$; in this case
$N_1(D_i)=\onespan{[W^i]}$ by corollary \ref{rho1}.\\
So if $\dim D_i=\dim \Di_{i}=3$ we have $\dim (D_i \cap \Di_{i}) \ge 1$, contradicting
the fact that $N_1(D_i)=\onespan{[W^i]}$ and $N_1(\Di_i)=\onespan{[\W^i]}$.
\end{proof}

\smallskip
\begin{lemma}\label{negative} 
Let $R_1$ be a divisorial ray of $X$, and $E_1$ its exceptional locus. 
If there exists a pair $(W^i,\W^{i}) \in \B$ such that $E_1 \cdot W^i <0$ and 
$E_1 \cdot \W^i>0$, then $[W^i] \in R_1$. 
\end{lemma} 

\begin{proof} 
Since $E_1 \cdot W^i <0$ we have $\loc(W^i) \subset E_1$; suppose 
by contradiction that $[W^i] \not \in R_1$.\\
If $\dim \loc(W^i_x) \ge 3$ for some $x$, then by lemma \ref{locy} (a) we have \linebreak 
$\dim \loc(W^i,R^1)_x \ge 5$, a contradiction since $W^i$ is noncovering.\\
It follows that $\dim \loc(W^i_x)=2$ and so, by lemma \ref{dimwi},
$\dim \loc(W^i)= 4$, hence $E_1 = \loc(W^i) = \loc(R^1,W^i)_x$ for some $x$; in particular 
by corollary \ref{ray}
\begin{equation*}
\cone(E_1)= \twospan{R_1}{[W^i]} \subset \cone(X)_{E_1<0}.
\end{equation*}
On the other hand, since $\dim \loc(\W^i_x) \ge 2$ and $E_1 \cdot \W^i >0$,
we have that $E_1$ contains curves proportional to $\W^i$, a contradiction. 
\end{proof}

\medskip
We can now resume the proof of theorem \ref{shelob}.\par
 
\medskip
{\bf Step 1.} \quad $\deg V = 4$.\par
\smallskip
Suppose by contradiction that $\deg V=5$, let $x \in X$ be a general point and
let $D$ be an irreducible component of $\loc(V_x)$; since $V$ is locally unsplit,
by corollary \ref{rho1} we have $N_1(D)=\onespan{[V]}$, and by proposition \ref{iowifam} 
we have $\dim D \ge \deg V -1 \ge 4$.\\
We are assuming $\rho_X \ge 2$, so it cannot be $D=X$, therefore $D$ is an effective divisor.\\
Thus the rc$\V$-fibration $\pi:\xymatrix@1{X~ \ar@{-->}[r]& ~Z}$ has fibers of dimension $\ge 4$; 
if $Z$ has positive dimension, take $V'$ to be a horizontal dominating family for $\pi$. 
By remark \ref{rcvdim} we know that
$\dim \loc(V'_x) = 1$, so $V'$ is covering and of degree $2$, hence it is unsplit, a contradiction.\\
So $X$ is rc$\V$-connected;
in particular $N_1(X)$ is generated as a vector space by the numerical
class of $V$ and the numerical classes of the families $W^i$ such that
$(W^i,\W^i) \in \B$ for some $\W^i$.\\
Consider the nonempty set of pairs $(W^i,\W^i) \in \B$ such that 
$[W^i] \neq [\alpha V]$, and the (non negative) intersection number $D \cdot V$.\\
If $D \cdot V=0$ then $D$ would be negative on at least one of these $W^i$
and so it would contain curves in $W^i$, against the fact that $N_1(D)=\onespan{[V]}$.\\
If else $D \cdot V >0$ then for every $i$ either $D \cdot W^i > 0$ or $D \cdot \W^i > 0$; but in 
this case either $\loc(W^i_x) \cap D$ or $\loc(\W^i_x) \cap D$ 
would be nonempty.\\
Suppose without loss of generality that we are in the first case; since, by lemma \ref{dimwi},
$\dim \loc(W^i_x) \ge 2$, then $\dim (\loc(W^i_x) \cap D)\ge 1$, against the fact 
that $N_1(\loc(W^i_x))=\onespan{[W^i]}$ and $N_1(D)=\onespan{[V]}$.\par

\smallskip
As a corollary of step 1 we get that $V$ is the unique locally unsplit dominating 
family for $X$ up to numerical equivalence: in fact, if $V'$ were another 
locally unsplit dominating family, for the general point $x \in X$ 
we would have $\dim (\loc(V_x) \cap \loc(V'_x)) \ge 1$
and so, since $N_1(\loc(V_x))=\onespan{[V]}$ and $N_1(\loc(V'_x))=\onespan{[V']}$, the families
would be proportional. But we have proved that $\deg V= \deg V'=4$, so $[V]=[V']$. \par

\medskip
{\bf Step 2}  \quad If $R_1$ is a small ray, then $[R^1]=[W^i]$ for some $i$. \par
\smallskip
By the results in Section 9, \cite{ACO}, we know that either $X$
is rc$\V$-connected or there exists an unsplit noncovering family
$V'$ such that $\deg V'=2$, $\dim \loc(V'_x)=2$ (and so, by proposition \ref{iowifam} 
$\dim \loc(V')=4$) and $X$ is rc$(\V,V')$-connected.\\
In particular, through every point of a fiber $F_1$ of $\f_{R_1}$
there passes a curve in $V$, in $W^i$ or in $V'$; if $V_y$ is unsplit for some $y \in F_1$
then $\dim(\loc(V_y) \cap F_1) \ge 1$, against the fact that
$N_1(V_y)=\onespan{[V]}$ and $N_1(F_1)=\onespan{R_1}$.\\
Therefore through every $y \in F_1$ there passes either a curve in $V'$ or a reducible cycle of $\V$; 
recalling that in $\B$ there is only a finite number of families we have that
either $F_1 \subset \loc(V')$ or  $F_1 \subset \loc(W^i)$ for some $i$.\\
Suppose first that $F_1 \subset \loc(V')$; note that, since $\loc(V')$ has dimension four
and $R_1$ is a small ray, the numerical class $[V']$ cannot belong to $R_1$,
hence $R^1$ is independent from $V'$ and we can apply lemma \ref{locy} 
(a) and obtain that $\dim \loc(V')_{F_1} \ge 5$, a contradiction.\\ 
If else $F_1 \subset \loc(W^i)$ and $R^1$ is independent from $W^i$ we reach
 a contradiction again by lemma \ref{locy} (a), so $[W^i]=[\alpha R^1]$.\\
By step 1 we know that $\deg W^i=2$, and since $R^1$ has minimal degree in the
ray we have also $\deg R^1 = 2$, i.e. $[R^1]=[W^i]$. \par

\medskip
{\bf Step 3} \quad $\rho_X=2$. \par
\smallskip
If $X$ is not rc{$\V$}-connected then the result follows from Section 9, Case 1 in \cite{ACO}.\\
Assume now that $X$ is rc{$\V$}-connected.\\
By the results in Section 9 of \cite{ACO}, if $\rho_X \ge 3$ then for every pair 
$(W^i,\W^{i}) \in \B$ and meeting components $D_i$ and $\Di_{i}$ of 
$\loc(W^i)$ and $\loc(\W^{i})$ we have 
$\dim D_i = \dim \Di_{i}=4$; moreover for at least one pair, say $(W^1,\W^1)$,
we know that $D_1 \neq \Di_1$.
In particular $X$ has no small contractions by step 2.\\
Let $E_1$ be the exceptional locus of a (divisorial) ray $R_1$ of $X$ and consider 
the intersection number $E_1 \cdot V$.\par

\smallskip
If $E_1 \cdot V>0$, for a general point $x \in X$
we have $\dim \loc(R^1)_{\loc(V_x)}=4$, so that
$E_1=\loc(R^1)_{\loc(V_x)}$ and $N_1(E_1)=\twospan{R_1}{[V]}$ by lemma \ref{numequns}.\\
If every pair $(W^i,\W^{i})$ lies in the plane spanned by $[V]$ and $[R^1]$ then
$\rho_X = 2$ and we conclude, otherwise let $(W^i, \W^{i})$ be a pair not lying in that plane:
then either $E_1 \cdot W^i >0$ or $E_1 \cdot \W^{i} >0$, implying that either
$E_1 \cap \loc(W^i_x)$ or $E_1 \cap \loc(\W^i_x)$ is nonempty and so has dimension $\ge 1$,
a contradiction.\par

\smallskip
If else  $E_j \cdot V=0$ for every $j$, then for every $j$
there exists  $i$ such that $E_j \cdot W^i < 0$ and $E_j \cdot \W^{i} >0$, 
so by lemma \ref{negative} $[W^i] \in R_j$ for some $i$ and $E_j = D_i$.\\
Let $R_k$ be an extremal ray such that $D_i \cdot R^k > 0$; by the argument above we know that
$E_k = \Di_l$ for some $l$, so it must be $D_i \cdot W^l <0$ and $i=l$, i.e. $E_k = \Di_i$.\\
It follows that $D_i \cdot V = \Di_i \cdot V=0$, but this is excluded in
Section 9, Case 2b of \cite{ACO} (note that since $D_i \cdot R^k >0$ we have
$D_i \neq E_k = \Di_i$).\par

\medskip
{\bf Step 4} \quad $X$ has a divisorial contraction.\par
\smallskip
Suppose that both the rays $R_1$ and $R_2$ correspond to small
contractions; by step 2 we know that there exist unsplit families
$W^1$ and $W^2$ such that $(R^1,W^1),(R^2,W^2) \in \B$, so in particular 
$[V]=[R^1]+[W^1] =[R^2]+[W^2]$.\\  
Let $x_1 \in \loc(R^1)$ and $x_2 \in \loc(R^2)$; since the contractions associated
to $R^1$ and $R^2$ are small, by inequality \ref{fiberlocus} we have
$\dim \loc(R^i)_{x_i} \ge 3$.\\
Denote by $D_i$ an irreducible component of $\loc(R^i,W^i)_{x_i}$;
since \\
$\loc(R^i,W^i)_{x_i}=\loc(W^i)_{\loc(R^i)_{x_i}}$, by 
lemma \ref{locy} (b) we have $\dim D_i=4$,
and by corollary \ref{ray} we have $\cone(D_1)=\twospan{R_1}{[W^1]}$, 
$\cone(D_2)=\twospan{R_2}{[W^2]}$.\\
It follows that $D_1 \cdot R^2 = D_2 \cdot R^1 =0$; moreover, since $D_i$ is an
effective divisor for every $i$, we have $D_i \cdot R^i >0$, so $D_i$ is nef.\\
Write $-K_X= a D_1 +bD_2$; we have
\begin{eqnarray*}
   4  & = & -K_X \cdot (W^1 +W^2)  =  \\
   & = & aD_1 \cdot (W^1 +W^2 +R^2) + bD_2 \cdot (W^1+W^2+R^1)= \\
   & = & aD_1 \cdot W^1 + aD_1 \cdot V + bD_2 \cdot W^2 + bD_2 \cdot V = \\
   & = & aD_1 \cdot W^1 + bD_2 \cdot W^2 -K_X \cdot V.
\end{eqnarray*}
Hence $aD_1 \cdot W^1 = bD_2 \cdot W^2 =0$, a contradiction. \par

\medskip
{\bf Step 5} \quad There exists a ray $R_1$ such that its associated contraction
$\f_1:X \to Y$ is a smooth blow-up of $Y$ along a smooth surface; moreover,
if $E_1$ is the exceptional divisor, $E_1 \cdot V >0$.
\smallskip

We know by step 4 that $X$ has a divisorial ray $R_1$; the other ray $R_2$ can be
either small or divisorial.\par

\smallskip
Let us start assuming that  $R_2$ is small; denote by $E_1$ the exceptional locus of 
$R_1$ and by $G_2$ a component of the exceptional locus of $R_2$.\\
The divisor $E_1$ is positive on $R^2$; it follows that $\loc(R^1)_{G_2}$ is nonempty
and has dimension four, so that $E_1 = \loc(R^1)_{G_2}$; in particular every
fiber of $R_1$ meets $G_2$ and so it is two-dimensional.\\
We can thus apply Theorem 5.1 of \cite{AO} and we get that $\f_1: X \to Y$ is a blow-down 
with center a smooth surface $S$ . \par

\smallskip
By step 2 there exists a pair $(W^1,\W^1) \in \B$ such that $[\W^1]=[R^2]$.
Take $D_1$ to be an irreducible
component of $\loc(W^1)$ which intersects $\loc(\W^1)$; since $D_1 \cdot \W^1 >0$ 
we have that $D_1 =\loc(\W^1,W^1)_x$ and so $\cone(D_1)=\twospan{[W^1]}{[\W^1]}$.\\
We claim that $[W^1] = [R^1]$: if this is not the case then $D_1 \neq E_1$, so 
$\f_1(D_1)$ is an effective divisor on $Y$. Moreover $\f_1(D_1)$ is ample since $\rho_Y=1$,
hence it meets $S$ and  $D_1 \cap E_1 \neq \emptyset$.\\ 
It follows that $\dim(D_1 \cap \loc(R^1_x)) \ge 1$ and
$D_1$ contains curves numerically proportional to $R^1$, a contradiction.\\
So we have proved that $\cone(X)=\twospan{R_1}{R_2}=\twospan{[W^1]}{[\W^1]}$;
moreover we have that $\loc(W^1)=E_1$.\par

\smallskip
Now we show that  $E_1 \cdot V >0$.\\
Suppose first that $X$ is not rc{$\V$}-connected and denote by $\pi$
the rc{$\V$}-fibration.
Let $X^0$ be the open subset on which $\pi$ is defined, take $x \in E_1 \cap X^0$ 
and consider $\loc(R^1_x)$; since $\dim Z \le 2$ either $\loc(R^1_x)$ is 
contained in a fiber of $\pi$ or $\loc(R^1_x)$ dominates $Z$.\\
In the first case, if $H$ is an ample divisor on $Z$ we have that 
$(\pi^*H) \cdot V = (\pi^*H) \cdot R^1 =0$,
so $\pi^*H$ is numerically trivial on $X$, a contradiction, hence 
$\loc(R^1_x)$ dominates $Z$ and $\dim Z=2$.\\
For a general $x$ in $X$ the fiber of $\pi$ through $x$ has dimension three and contains 
$\loc(V_x)$, hence $F_x=\loc(V_x)$; $E_1$ meets this fiber and cannot contain it,
so $E_1 \cdot V>0$.\par

\medskip
Assume now that $X$ is rc{$\V$}-connected and
suppose by contradiction that $E_1 \cdot V = 0$; in this case, by lemma \ref{negative}
$\B$ contains only the pair $(W^1,\W^1)$ and possibly pairs
$(W^j,\W^j)$ with $[W^j] =[\W^j] = \frac{1}{2}[V]$.\\
Let $T = \loc(W^1) \cup \loc(\W^1)$ and take a point $x$ outside $T$; since
$X$ is rc$\V$-connected we can join $x$ and $T$ with a chain of cycles in $\V$.
Let $\Gamma$ be the first irreducible component which meets $T$.\\
Clearly $\Gamma$ cannot belong to $[W^1]$ and $[\W^1]$ because it is not contained in $T$,
so it belongs either to $V$ or to $W^j$ for some $j$, say $j=2$. Since $E_1 \cdot V = E_1 \cdot W^2=0$,
$\Gamma$ must intersect $T$ in points of $T \setminus \loc(W^1)$.\\
Let $y$ be a point in $\Gamma \cap T$ and let $G_y$ be the irreducible 
component of $T$ which contains $y$; by Lemma 9.1 in \cite{ACO} we have 
that either $\Gamma \subset \loc(V_z)$
for some $z$ such that $V_z$ is unsplit or $\Gamma \subset \loc(W^2)$.\\ 
In the first case we have $\dim (\loc(V_z) \cap G_y) \ge 1$, against the fact that 
$N_1(V_z)=\onespan{[V]}$ and $N_1(G_y)=\onespan{[\W^1]}$.\\
In the second case we consider $H = \loc(W^2)_{G_y}$: by lemma \ref{locy} we have
$\dim H = 4$, and by lemma \ref{ray} we have $\cone(H) = \twospan{[\W^1]}{[W^2]}$; this implies that
$H \cdot R^1 = 0$.\\
The image $\f_1(H)$ of $H$ in $Y$ is an effective, hence ample, divisor;
therefore $\f_1(H) \cap S \not = \emptyset$ and $H \cap E_1 \neq \emptyset$.\\
For every point $t \in H \cap E_1$ we have that both $\loc(W^2_t)$ and $\loc(W^1_t)$ are
contained in $H \cap E_1$, since $H \cdot W^1 = E_1 \cdot W^2=0$.\\
Therefore $\loc(W^2,W^1)_t \subseteq H \cap E_1$, 
and we reach a contradiction since $\dim \loc(W^2,W^1)_t = 4$ by lemma \ref{locy} (a).\par

\medskip
Assume now that both $R_1$ and $R_2$ are divisorial and let $E_1$, $E_2$ be the respective
exceptional loci.\\
We cannot have $E_1 \cdot V = E_2 \cdot V=0$ (see the end of proof of Theorem 7.1 in \cite{ACO}),
so we can suppose that $E_1 \cdot V > 0$.\\
If $x \in X$ is a general point then $\loc(R^1)_{\loc(V_x)}$ is nonempty
and has dimension four, so $E_1 = \loc(R^1)_{\loc(V_x)}$; in particular every
fiber of $R_1$ meets $\loc(V_x)$ and so it is two-dimensional.\\
Now we apply Theorem 5.1 of \cite{AO} and we get that $\f_1: X \to Y$ is a blow-down 
with center a smooth surface $S$.\par

\medskip
{\bf Step 6} \quad $Y \simeq \proj^5$.\par
\smallskip

Let $V_Y$ be a minimal covering family for $Y$ and let $V^*_Y$ be the
family of deformations of the strict transform of a general curve in $V_Y$. We know that 
$V^*_Y$ is covering and that
\begin{equation} \label{blowup}
6 = \dim Y+1 \ge -\f_1^*K_Y \cdot V_Y^* = -K_X \cdot V_Y^* +2E_1 \cdot V_Y^*.
\end{equation}
By Proposition 3.7 in \cite{Kob} a general member of $V_Y$ does not meet $S$, which
has codimension three in $Y$, hence $E_1 \cdot V^*_Y =0$.\\
The family $V^*_Y$ cannot be locally unsplit: otherwise, by the final part of step one,
we would have $[V^*_Y] = [V]$, but we know by step 5 that $E_1 \cdot V>0$.\\
It follows that one among the families of
irreducible components of cycles in $\V_Y^*$ is covering, and so, by our assumptions,
the degree of this family is at least four; in particular $-K_X \cdot V^*_Y \ge 6$.
Equation \ref{blowup} yields
$$6 \ge  -K_Y \cdot V_Y -\f_1^*K_Y \cdot V_Y^* = -K_X \cdot V_Y^* \ge 6,$$  
so we can conclude that $Y$ has a minimal
(hence locally unsplit by lemma \ref{horizontal}) covering family of degree
$6=\dim Y+1$; by the proof of Theorem 1.1 in \cite{Ke} we have $Y \simeq \proj^5$. (Note
that the assumptions of the quoted result are different, but
the proof actually works in our case since for a very general $y$ the pointed family 
$(V_Y)_y$ has the properties 1-3 in Theorem 2.1 of \cite{Ke}).\par

\medskip
{\bf Final step}\par
\smallskip

Let $S \subset \proj^5$ be the center of the blow-up, let $l$ be a (bi)secant line of $S$
and let $\tilde l$ be the proper transform of $l$; then
\begin{equation*}
-K_X \cdot \tilde l = \f_1^*\Ol(6) \cdot \tilde l -2E_1 \cdot \tilde l =2,
\end{equation*}
so the corresponding family on $X$ is unsplit. Since $X$ does not admit unsplit covering 
families, through the general point of $\proj^5$ there is no secant line of $S$.\\
It follows that either $S$ is a Veronese surface or $S$ is degenerate.\\
If $S$ is contained in an hyperplane $H$ and in no three-dimensional linear
subspace of $\proj^5$, through every point of $H$ there is 
a secant line, so the strict transform $\acca$ of $H$ is the locus of an unsplit family $W$
on $X$.\\
Since $\f_1^*\Ol(1) = \f_1^* H = \acca + k E_1$ with $k >0$ we have
$\acca = \f_1^* \Ol(1) -kE_1$, so $\acca \cdot W <0$. 
It follows that $\acca $ is negative on $R_2$ and it follows that $\acca =E_2$, $R_2$ is divisorial and
$W=R^2$ by lemma \ref{negative}.\\
Again from the canonical bundle formula we have $E_1 \cdot R^2=2$, so $k = 1$ and 
$E_2 = \acca = \f_1^*\Ol(1) - E_1$; in particular $\f_1^*\Ol(1) \cdot R^2 = 1$ and
the contraction $\f_2: X \to Z$ is supported by $K_X + 2\f_1^*\Ol(1) + \f_2^*A $ for some very ample divisor 
$A$ on $Z$. 
By Corollary 5.8.1 of \cite{AW} $\f_2$ is equidimensional, so it is a smooth blow-down.\\
Computing the canonical bundle of $E_2$ 
$${K}_{E_2} = -5 \f_1^* \Ol(1) +E_1,$$
we find that $E_2$ is a Fano variety.\\
Moreover $E_2$ has a $\proj^2$-bundle structure over a smooth surface $S' \subset Z$, and
since $\rho(E_2)=2$ we have $S' \simeq \proj^2$; by the classification in
\cite{SW} we know that $S$ is a cubic scroll.\\
Finally, if $S$ is contained in a three-dimensional linear subspace $\Lambda \subset \proj^5$ and
$l$ is a line in $\Lambda$ we have
\begin{equation*}
-K_X \cdot \tilde l = \f_1^*\Ol(6) \cdot \tilde l -2E_1 \cdot \tilde l = 6 - 2\deg S,
\end{equation*}
so $S$ has degree $\le 2$ and it cannot be
a plane, since the blow-up of $\proj^5$ along a two-dimensional plane
has a fiber type contraction.
\end{proof}


\section{The higher-dimensional case}

Let $V$ be a covering locally unsplit family of rational curves on a Fano variety $X$ 
of dimension $n \ge 6$ and 
pseudoindex $i_X = n - 3$; 
as shown in the case of fivefolds, if $\rho_X \ge 2$ then $\deg V \le n$.\\

\begin{lemma}
Let $X$ be a Fano variety of dimension $n \ge 6$, pseudoindex $i_X = n - 3$ and 
Picard number $\rho_X \ge 2$; let $V$ be a covering locally unsplit family of rational curves on $X$. 
Then $V$ is unsplit.
\end{lemma}

\begin{proof}
If $\deg V \le n - 1 < 2i_X$ then $V$ is unsplit, so
we can assume that $\deg V = n$.
Let $x \in X$ be a general point and let $D$ be an irreducible component of $\loc(V_x)$; 
since $V$ is locally unsplit we have $N_1(D)=\onespan{[V]}$ by corollary \ref{rho1} and,
by proposition \ref{iowifam}, $\dim D \ge \deg V -1 = n - 1$.
We are assuming $\rho_X \ge 2$, so it cannot be $D = X$, therefore $D$ is an effective divisor.\\
The rc$\V$-fibration $\pi:\xymatrix@1{X~ \ar@{-->}[r]& ~Z}$ has fibers of dimension $\ge n - 1$; if $Z$ has positive
dimension, take $V'$ to be a horizontal dominating family.
Then if $F$ is a fiber of $\pi$ we have
$$\dim(F \cap \loc(V'_x)) \ge n - 1 + \deg V' - 1 - n \ge i_X - 2 \ge 1,$$
contradicting lemma \ref{horizontal}.
So $X$ is rc$\V$-connected, and we reach a contradiction as in step 1 of section \ref{without}. 
\end{proof}

\smallskip
\begin{corollary}
In the above assumptions, Mukai conjecture (\ref{conj}) holds for $X$. 
In particular if $\dim X \ge 8$ then $\rho_X = 1$ except if $X \simeq \proj^4 \times \proj^4$.
\end{corollary}

\begin{proof}
Since $i_X \ge \frac{n+3}{3}$ we can apply Theorem 7.1 in \cite{ACO}.
\end{proof}

\medskip
To conclude the proof of theorem \ref{main} we have to deal with varieties
of dimension $6$ and $7$.
Note that Mukai conjecture implies that $\rho_X \le 2$ except if 
$X \simeq \pdue \times \pdue \times \pdue$.\par

\medskip
Arguing as in case $\rho_X = 2$ of section \ref{with} we can prove that
$X$ has at least one fiber type ray $R_1$. 
Let $R_2$ be the other extremal ray; then if $F_i$ is a general fiber of
the contraction $\f_{R_i}$ we know that
   $$\dim F_1 + \dim F_2 \le \dim X,$$
and together with the fiber locus inequality this concludes the classification.

\section*{Acknowledgements}

We thank Marco Andreatta, Edoardo Ballico, Massimiliano Mella and 
Roberto Pi\-gna\-telli for helpful remarks and discussions.
We also thank the referee for pointing out some inaccuracies in the original version.


\bibliographystyle{amsplain}

\end{document}